\documentclass[a4paper,12pt]{amsart}
\usepackage[utf8]{inputenc}
\usepackage{amsmath,amssymb, comment,hyperref}
\usepackage{multirow,array,enumerate,wasysym}
\usepackage[left=2.9cm,right=2.9cm,bottom=2.9cm,top=2.9cm]{geometry}
\usepackage{amsthm}
\usepackage{natbib}
\usepackage{tikz}
\usepackage{microtype}

\bibpunct{[}{]}{,}{n}{}{;}

\newtheorem{theorem}{Theorem}[section]
\newtheorem{lemma}[theorem]{Lemma}

\newtheorem{corollary}[theorem]{Corollary}
\theoremstyle{definition}

\newtheorem{example}[theorem]{Example}
\newtheorem{remark}[theorem]{Remark}

\newcommand{\SZ}{\mathbb{Z}}                    
\newcommand{\SC}{\mathbb{C}}                    
\newcommand{\CF}{\mathcal{F}}                    

\newcommand{\ra}[1]{\kern-1.5ex\xrightarrow{\ \ #1\ \ }\phantom{}\kern-1.5ex}
\newcommand{\ras}[1]{\kern-1.5ex\xrightarrow{\ \ \smash{#1}\ \ }\phantom{}\kern-1.5ex}

\title{Hilbert scheme of points on cyclic quotient singularities of type $(p,1)$}
\author{Ádám Gyenge}

\begin{document}
\begin{abstract}
 In this note we investigate the generating series of the Euler characteristics of Hilbert scheme of points on cyclic quotient singularities of type $(p,1)$. We link the appearing combinatorics to $p$-fountains, a generalization of the notion of fountain of coins. We obtain a representation of the generating series as coefficient of a two variable generating series.
\end{abstract}

\maketitle


\section{Introduction}
\subsection{Hilbert schemes of points}
For a complex quasi-projective surface $X$ and a positive integer $n$, let $\mathrm{Hilb}^n(X)$ be the Hilbert scheme of $n$ points on $X$ parametrizing the zero-dimensional subschemes of $X$ of length $n$. As a set, it consists of ideals $I$ of the ring $\mathcal{O}_X$ of regular functions on $X$, such that $\mathcal{O}_X/I$ is an $n$-dimensional $\SC$-vector space. 
The disjoint union of these spaces for all $n$ is denoted as $\mathrm{Hilb}(X)$. 
In this paper we are interested in the topological Euler-characteristics of the Hilbert schemes of points on $X$. We collect these in a generating series:
\[Z_{X}(q)=\sum_{n\geq 0} q^n\chi(\mathrm{Hilb}^n(X)). \]

The Hilbert scheme of points on smooth surfaces is an object investigated for a long time. It has turned out that it carries a rich geometrical structure, see e.g. \cite{nakajima1999lectures}. By a theorem of Fogarty, it is a smooth variety. 
The following result of Göttsche gives an expression of the generating series of the Poincaré polynomial of the Hilbert schemes of points on $X$ in terms of the Betti numbers of $X$.
\begin{theorem}[\cite{gottsche1990betti}] If $X$ is a smooth quasi-projective surface, then
\[\sum_{l=1}^\infty q^l P_t( X^{[l]} )=\prod_{m=1}^\infty \frac{(1+t^{2m-1}q^m)^{b_1(X)}(1-t^{2m+1}q^m)^{b_3(X)}}{(1-t^{2m-2}q^m)^{b_0(X)}(1-t^{2m}q^m)^{b_2(X)}(1-t^{2m+2}q^m)^{b_4(X)}} ,\]
where $b_i(X)$ is the $i$-th Betti number of $X$, and $P_t(Y)=\sum_{i \geq 0}t^ib_i(Y)$. 
\end{theorem}
Applying this to $X=\SC^2$ and substituting $t=-1$ we obtain
\begin{corollary}
\label{cor:planepart}
\[Z_{\SC^2}(q)=\prod_{m=1}^\infty \frac{1}{1-q^m}. \]
\end{corollary}

\subsection{Cyclic quotient singularities and orbifolds}

Fix a positive integer $p$. Let $\SZ_p$ be the cyclic group of order $p$ with generator $g$ and let it act on $\SC^2$ as:  $g.x=\mathrm{e}^{\frac{2 \pi i}{p}} x$, $g.y=\mathrm{e}^{\frac{2 \pi i q}{p}}y$ where $q$ is coprime to $p$. Then we get an action of $\SZ_p$ on $\SC^2$ which is totally discontinuous away from the origin. Let $X(p,q)$ denote the quotient variety. It is called the cyclic quotient singularity of type $(p,q)$.

The Hilbert scheme $\mathrm{Hilb}(X(p,q))$ of points on $X(p,q)$  is the moduli space of ideals in $\mathcal{O}_{X(p,q)}=\SC[x,y]^{\SZ_p}$ of finite colength. The case when $q=p-1$ is called the type $A$ singularity. It was considered in \cite{gyenge2015euler}.
We are now interested in the other extreme case when $q=1$, that is, in calculating
\[Z_{X(p,1)}(q)=\sum_{m=0}^\infty \chi\left(\mathrm{Hilb}^m(X(p,1))\right)q^m.\]



\subsection{The result}

Our main result is a 
representation of $Z_{X(p,1)}(q)$ as coefficient of a two variable generating function. In these two variable generating function continued fractions appear. We introduce the notation $[z^0]\sum_n A_n z^n=A_0$.
\begin{theorem} 
\label{thm:main}
Let $F(q,z)$ and $T(q,z)$ be the functions defined in \eqref{eq:Fdef} and \eqref{eq:Tdef} below, respectively. Then:
\[ Z_{X(p,1)}(q)=[z^0]T(q,z)\left(F(q^{-1},z^{-1})-(qz)^{-p}F(q^{-1},(qz)^{-1})\right).\]
\end{theorem}

\begin{remark}
\begin{enumerate}
\item  For $p=1$, $X(1,1)=\SC^2$. Then, by Corollary \ref{cor:planepart}, $Z_{X(1,1)}(q)$ is also equal to 
\[\prod_{m=1}^\infty \frac{1}{1-q^m}.\]
\item For $p=2$, $Z_{X(2,1)}(q)$ was calculated in \cite{gyenge2015euler} in a different way. It follows from those results, that $Z_{X(2,1)}(q)$ is also equal to 
\[ \left( \prod_{m=1}^\infty (1-q^m)^{-1} \right)^{2} \cdot\sum_{ m\in \SZ } \xi^{m}q^{m^2}\;,\]
where $\xi=\mathrm{exp}(\frac{2 \pi i}{3})$.
\end{enumerate}
\end{remark}

\section{Proof}

\subsection{Torus action on \texorpdfstring{$X(p,1)$}{X(p,1)}}

Let us fix $p$ once and for all. The surface singularity $X(p,1)$ is toric, i.e. it carries a $(\SC^{\ast})^2$-action with an isolated fixed point. On the level of regular functions the fixed points are the monomials in
\[A:=\mathcal{O}_{X(p,1)}=\SC[x,y]^{\SZ_p}\cong\SC[x^p,x^{p-1}y,\dots,xy^{p-1},y^{p}].\] The $(\SC^\ast)^2$-action lifts to $\mathrm{Hilb}^n(X(p,1))$ for each $n$. The action of $(\SC^\ast)^2$ on $\mathrm{Hilb}^n(X(p,1))$ again has only isolated fixed points which are given by the finite colength monomial ideals in $A$.

It is well known that finite colength monomial ideals inside $\SC[x,y]$ are in one-to-one correspondence with partitions and with Young diagrams. The generators of a monomial ideal are the functions corresponding to those blocks in the complement of the diagram which are at the corners. Since $A\subset \SC[x,y]$, each fixed point of the $(\SC^\ast)^2$-action on $\mathrm{Hilb}^n(X(p,1))$ corresponds to a Young diagram.

A block at position $(i,j)$ is called a $0$-block if $i+j\equiv 0\; (\textrm{mod}\; n)$. We will call a Young diagram $0$-generated if its generator blocks are $0$-blocks. The $0$-weight of a (not necessarily $0$-generated) Young diagram $\lambda$ is the number of $0$-blocks inside the Young diagram. It is denoted as $\mathrm{wt}_0(\lambda)$. Let $\mathcal{P}_0$ be the set of $0$-generated Young diagrams. It decomposes as
\[\mathcal{P}_0=\bigsqcup_{n\geq 0} \mathcal{P}_0(n), \]
where $\mathcal{P}_0(n)$ is the set of $0$-generated Young diagrams which have $0$-weight $n$. 

\begin{example} \label{ex:1}
If $p=3$, then the following is a 0-generated Young diagram (we also indicated the generating 0-blocks):

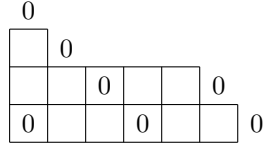
\begin{center}
\begin{tikzpicture}[scale=0.5, font=\footnotesize, fill=black!20]
\draw (0, 0) -- (6,0);
\draw (0,1) --(6,1);
\draw (0,2) --(5,2);
\draw (0,3) --(1,3);
\draw (0,0) -- (0,3);
\draw (1,0) -- (1,3);
\draw (2,0) -- (2,2);
\draw (3,0) -- (3,2);
\draw (4,0) -- (4,2);
\draw (5,0) -- (5,2);
\draw (6,0) -- (6,1);
\draw (0.5,0.5) node {0};
\draw (3.5,0.5) node {0};
\draw (2.5,1.5) node {0};

\draw (1.5,2.5) node {0};
\draw (0.5,3.5) node {0};
\draw (6.5,0.5) node {0};
\draw (5.5,1.5) node {0};
\end{tikzpicture}
\end{center}
The 0-weight of this Young diagram is 3.
\end{example}

The generators of any ideal of $A$, when considered as functions in $\SC[x,y]$, have to be invariant under the $\SZ_p$-action. This implies the following statement.
\begin{lemma}The monomial ideals inside $A$ 
of colength $n$ are in one-to-one correspondence with the $0$-generated Young diagrams in $\mathcal{P}_0(n)$.
\end{lemma}
\begin{corollary}
\label{cor:sumdiag}
\[Z_{X(p,1)}(q)=\sum_{\lambda \in \mathcal{P}_0}q^{\mathrm{wt}_0(\lambda)}=\sum_{n\geq 0} |\mathcal{P}_0(n)|q^n,\]
where $|S|$ denotes the number of elements in the set $S$.
\end{corollary}

\subsection{Augmentation of 0-generated Young diagrams}
\label{subsec:aug}


For any 0-generated Young diagram, we can consider the area which is between the diagram and the line $x+y=c$. Here $c$ is an integer congruent to 0 modulo $p$ and it is as small as possible such that the line $x+y=c$ does not intersect the diagram. In other words, this line is just the antidiagonal closest to the diagram.  The line cuts out the smallest isosceles right-angled triangle which contains the whole diagram. In the case of Example \ref{ex:1}, this looks as follows:
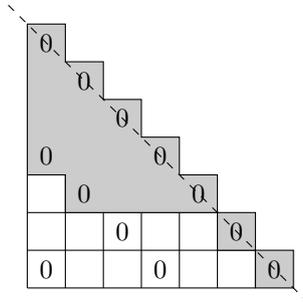
\begin{center}
\begin{tikzpicture}[scale=0.5, font=\footnotesize, fill=black!20]
\draw (0, 0) -- (6,0);
\draw (0,1) --(6,1);
\draw (0,2) --(5,2);
\draw (0,3) --(1,3);
\draw (0,0) -- (0,3);
\draw (1,0) -- (1,3);
\draw (2,0) -- (2,2);
\draw (3,0) -- (3,2);
\draw (4,0) -- (4,2);
\draw (5,0) -- (5,2);
\draw (6,0) -- (6,1);
\draw (0.5,0.5) node {0};
\draw (3.5,0.5) node {0};
\draw (2.5,1.5) node {0};

\filldraw (6,0) -- (6,1) -- (5,1) -- (5,2) -- (1,2) -- (1,3) -- (0,3) -- (0,7) -- (1,7) -- (1,6) -- (2,6) -- (2,5) -- (3,5) -- (3,4) -- (4,4) -- (4,3) -- (5,3) -- (5,2) -- (6,2) -- (6,1) -- (7,1) -- (7,0)--(6,0) --cycle  ;
\draw (1.5,2.5) node {0};
\draw (0.5,3.5) node {0};
\draw (6.5,0.5) node {0};
\draw (5.5,1.5) node {0};
\draw (4.5,2.5) node {0};
\draw (3.5,3.5) node {0};
\draw (2.5,4.5) node {0};
\draw (1.5,5.5) node {0};
\draw (0.5,6.5) node {0};

\draw[dashed] (-0.5,7.5) -- (7.5,-0.5);

\end{tikzpicture}
\end{center}

If, for example, $p=1$ then the area between the diagram, the $x$ and $y$ coordinate axes, and the above mentioned antidiagonal, when rotated 45 degrees counterclockwise and flipped, is a special type of a \textit{fountain of coins} as introduced in \cite{odlyzko1988editor}. An \emph{$(n,k)$ fountain (of coins)} is an arrangement of $n$ coins in rows such that there are exactly $k$ consecutive coins in the bottom row, and such that each coin in a higher row touches exactly two coins in the next lower row. In our case, the $0$-blocks in the area under consideration (which is colored grey in the diagram above) are replaced by coins or circles or zero symbols.

We generalize this notion to arbitrary $p$. An \emph{$(n,k)$ $p$-fountain} is an arrangement of $n$ coins in rows such that 
\begin{itemize}
\item there are exactly $k$ consecutive coins in the bottom row,
\item immediately below each coin there are exactly $p$+1 ``descendant'' coins in the next lower row,
\item and $p$ coins among the descendants of two neighboring coins coincide.
\end{itemize}
In other words, we rotate and flip the area between the axes and a specific antidiagonal, and the coins can be placed exactly on the 0-blocks, such that if there is a coin somewhere, then there has to be coins on all the 0-blocks in the area which have higher $x$ or $y$ coordinates in the original orientation of the plane. The empty diagram is considered as a $(0,0)$ $p$-fountain.

Continuing the case of Example \ref{ex:1} further, the associated $(9,7)$ $3$-fountain is
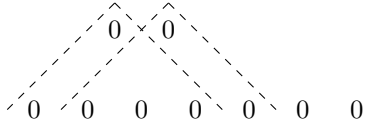
\begin{center}
\begin{tikzpicture}[scale=0.5, font=\footnotesize, fill=black!20, rotate=-45, xscale=1,yscale=-1]

\draw[dashed] (0,3) -- (4,3);
\draw[dashed] (0,3) -- (0,7);

\draw[dashed] (1,2) -- (5,2);
\draw[dashed] (1,2) -- (1,6);

\draw (1.5,2.5) node {0};
\draw (0.5,3.5) node {0};
\draw (6.5,0.5) node {0};
\draw (5.5,1.5) node {0};
\draw (4.5,2.5) node {0};
\draw (3.5,3.5) node {0};
\draw (2.5,4.5) node {0};
\draw (1.5,5.5) node {0};
\draw (0.5,6.5) node {0};


\end{tikzpicture}
\end{center}
Here we also indicated the descendants of the coins in the upper row.

An $(n,k)$ $p$-fountain  is called primitive if its next-to-bottom row contains no empty positions, i.e. contains $k-p$ coins. In particular, the fountain with $p$ coins in the bottom row but with no coin in the higher rows is primitive, but the ones with less than $p$ coins in the bottom row are not primitive. The fountains that appear between our 0-generated Young diagrams and the diagonals are special because in each case there is at least one empty position in the next-to-bottom row, so they correspond exactly to the non-primitive $p$-fountains.  

Let $f(n,k)$ (reps., $g(n,k)$) be the number of arbitrary (resp., primitive )$(n,k)$ $p$-fountains. Let $F(q,z)=\sum_{n,k\geq 0}f(n,k)q^nz^k$ (resp., $G(q,z)=\sum_{n,k\geq 0}g(n,k)q^nz^k$) be the two variable generating function of the sequence $f(n,k)$ (resp., $g(n,k)$). We will calculate $F(q,z)$ and $G(q,z)$ by extending the ideas of \cite{odlyzko1988editor}.

By removing the bottom row of a primitive $(n,k)$ $p$-fountain one obtains a $(n-k,k-p)$ $p$-fountain. Therefore,
\[ g(n,k)=f(n-k,k-p) \qquad (n\geq k, k\geq p), \]
and
\begin{equation}\label{eq:fshifted} G(q,z)=(qz)^{p}F(q,qz).\end{equation}
We prescribe that 
\begin{equation} \label{eq:initcond}f(0,0)=\dots=f(p-1,p-1)=1 \textrm{ and } g(0,0)=\dots=g(p-1,p-1)=0.\end{equation}

Let us consider an arbitrary $(n,k)$ $p$-fountain $\CF$, and assume that the first empty position in the next-to-bottom row is the $r$-th ($1 \leq r \leq k-p+1$). Then we can split $\CF$ into a primitive $(m+p-1,r+p-1)$ $p$-fountain and a not necessarily primitive $(n-m,k-r)$ $p$-fountain after the first and before the last descendant coin of the above mentioned missing $r$-th position. The descendant coins at the second to the penultimate positions will be doubled. 

Sticking to our favorite Example \ref{ex:1}, the splitting looks as follows:

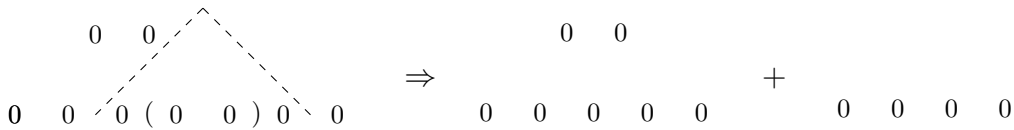
\begin{center}
\begin{tabular}{m{5cm} m{0.5cm} m{3.5cm} m{0.5cm} m{3.5cm} }
\begin{tikzpicture}[scale=0.5, font=\footnotesize, fill=black!20, rotate=-45, xscale=1,yscale=-1]
\draw[dashed] (2,1) -- (6,1);
\draw[dashed] (2,1) -- (2,5);
\draw (1.5,2.5) node {0};
\draw (0.5,3.5) node {0};
\draw (6.5,0.5) node {0};
\draw (5.5,1.5) node {0};
\draw (4.5,2.5) node {0};
\draw (3.5,3.5) node {0};
\draw (2.5,4.5) node {0};
\draw (1.5,5.5) node {0};
\draw (0.5,6.5) node {0};
\draw (0.5,6.5) node {0};
\draw (3,4) node {(};
\draw (5,2) node {)};
\end{tikzpicture}
&
$\Rightarrow$
&
\begin{tikzpicture}[scale=0.5, font=\footnotesize, fill=black!20, rotate=-45, xscale=1,yscale=-1]
\draw (1.5,2.5) node {0};
\draw (0.5,3.5) node {0};
\draw (4.5,2.5) node {0};
\draw (3.5,3.5) node {0};
\draw (2.5,4.5) node {0};
\draw (1.5,5.5) node {0};
\draw (0.5,6.5) node {0};
\end{tikzpicture}
&+
&
\begin{tikzpicture}[scale=0.5, font=\footnotesize, fill=black!20, rotate=-45, xscale=1,yscale=-1]
\draw (1.5,2.5) node {};
\draw (3.5,3.5) node {0};
\draw (2.5,4.5) node {0};
\draw (1.5,5.5) node {0};
\draw (0.5,6.5) node {0};
\end{tikzpicture}
\end{tabular}
\end{center}
The dashed line on the left picture indicates the missing coin from the second row and its descendants. The primitive part is to the left of the right parenthesis, while the remaining part is to the right of the left parenthesis.

This factorization is unique. 
\begin{equation} \label{eq:conv} f(n,k)=\sum_{\substack{0\leq m \leq n-p+1 \\ 0 \leq r \leq k-p+1}} g(m+p-1,r+p-1)f(n-m,k-r) \qquad (n,k\geq p). \end{equation}
 The conditions $m \leq n-p+1$ and $r \leq k-p+1$ are equivalent to $p-1 \leq n-m$ and $p-1 \leq k-r$ respectively. That is, the other remaining fountain has to have a first row of length at least $p-1$. From \eqref{eq:initcond} we see that $F(q,z)$ has the form
 \[F(q,z) =1+qz+\dots+ (qz)^{p-1}+\dots.\]
 Then $F(q,z)-1-qz-\dots-(qz)^{p-2}$ is exactly the generating function of fountains, which have a first row of length at least $p-1$, that is, which can appear as the second factor on the right hand side of \eqref{eq:conv}. Therefore,
 \[(qz)^{-p+1}G(q,z)(F(q,z)-1-qz-\dots-(qz)^{p-2})\]
 enumerates all $p$-fountains for which $n,k \geq p$. The generating function of all $p$-fountains then satisfies
\begin{equation} \label{eq:frec} 
F(q,z) =1+qz+\dots+ (qz)^{p-1}+(qz)^{-p+1}G(q,z)(F(q,z)-1-qz-\dots-(qz)^{p-2}). 
\end{equation}

Using \eqref{eq:frec},
\[ (F(q,z)-1-qz-\dots-(qz)^{p-2})(1-(qz)^{-p+1}G(q,z))=(qz)^{p-1}.  \]
Then, by \eqref{eq:fshifted}, the generating function of $p$-fountains is
\begin{equation}\label{eq:Fdef} \begin{aligned}F(q,z)&=\frac{(qz)^{p-1}}{1-(qz)^{-p+1}G(q,z)}+1+qz+\dots+(qz)^{p-2}\\&=\frac{(qz)^{p-1}}{1-qzF(q,qz)}+\frac{1-(qz)^{p-1}}{1-qz}\\
&=\frac{(qz)^{p-1}}{1-qz\left(\frac{(q^2z)^{p-1}}{1-q^2zF(q,q^2z)}+\frac{1-(q^2z)^{p-1}}{1-q^2z}\right)}+\frac{1-(qz)^{p-1}}{1-qz}=\dots\\&=\frac{(qz)^{p-1}}{1-qz\left(\frac{(q^2z)^{p-1}}{1-q^2z\left(\frac{(q^3z)^{p-1}}{1-\dots}+\frac{1-(q^3z)^{p-1}}{1-q^3z} \right)}+\frac{1-(q^2z)^{p-1}}{1-q^2z}\right)}+\frac{1-(qz)^{p-1}}{1-qz}.
\end{aligned}  \end{equation}
Consequently, the generating function of primitive $p$-fountains is
\begin{equation}\label{eq:Gdef} G(q,z)=\frac{(q^2z)^{p-1}}{1-q^2z\left(\frac{(q^3z)^{p-1}}{1-q^3z\left(\frac{(q^4z)^{p-1}}{1-\dots}+\frac{1-(q^4z)^{p-1}}{1-q^4z} \right)}+\frac{1-(q^3z)^{p-1}}{1-q^3z}\right)}+\frac{1-(q^2z)^{p-1}}{1-q^2z}.\end{equation}

\begin{remark}
For $p=1$, Ramanujan \cite[p. 104]{andrews1998theory} obtained the beautiful formula
\[ F(q,z)= \frac{1}{1-\frac{qz}{1-\frac{q^2z}{1-\dots}}}=\frac{\sum_{n\geq 0} (-qz)^n \frac{q^{n^2}}{(1-q)(1-q^2)\dots(1-q^n)}}{\sum_{n\geq 0} (-z)^n \frac{q^{n^2}}{(1-q)(1-q^2)\dots(1-q^n)}}.\]
\end{remark}

The number of non-primitive $(n,k)$ $p$-fountains is obviously
\[ h(n,k)=f(n,k)-g(n,k),\]
which gives
\begin{equation}\label{eq:Hexpr1} H(q,z)=F(q,z)-G(q,z)=F(q,z)-(qz)^{p}F(q,qz) \end{equation}
for the generating series $H(q,z)$ of the numbers $h(n,k)$.


\subsection{Completion of the proof}

As mentioned in \ref{subsec:aug}, we augment each 0-generated Young diagram to the smallest  isosceles right-angled triangle. The area between the diagram and the triangle will be a $p$-fountain. For a fixed $p$, the hypotenuse of the possible isosceles right-angled triangles contains $lp+1$ blocks, where $l$ is a non-negative integer. The number of 0-blocks in the triangle with $lp+1$ blocks on the hypotenuse is $\sum_{0 \leq i \leq l}ip+1=p\frac{l(l+1)}{2}+l+1$. Therefore, the two variable generating series of these triangles is
\[\sum_{l\geq 0} q^{p\frac{l(l+1)}{2}+l+1}z^{lp+1}.\]
We will see immediately, that adding the terms with negative $l$ will not affect the final result. Hence, we define
\begin{equation}\label{eq:Tdef} 
\begin{aligned}
T(q,z)& =\sum_{l=-\infty}^{\infty} q^{p\frac{l(l+1)}{2}+l+1}z^{lp+1}=(qz)\sum_{l=-\infty}^{\infty}(q^p)^{\binom{l+1}{2}}(qz^p)^l\\
& =(qz)\prod_{n=1}^{\infty}(1+z^pq^{np+1})(1+z^{-p}q^{(n-1)p-1})(1-q^{np}),
\end{aligned}
\end{equation}
where at the last equality we have used the following form of the Jacobi triple product identity:
\[ \prod_{n=1}^{\infty}(1+zq^n)(1+z^{-1}q^{n-1})(1-q^n)=\sum_{j=-\infty}^{\infty}z^jq^{\binom{j+1}{2}}. \]

The generating series of 0-generated Young diagrams is then
\begin{equation}\label{eq:z0coeff} \sum_{\lambda \in \mathcal{P}_0}q^{\mathrm{wt}_0(\lambda)}=[z^0]T(q,z)H(q^{-1},z^{-1}).\end{equation}
Taking the coefficient of $z^0$ ensures that the hypotenuse of the triangle and the bottom row of the $p$-fountain match together, and also that the terms of $T(q,z)$ with negative powers of $z$ do not contribute into the result. Putting together Corollary \ref{cor:sumdiag}, \eqref{eq:Hexpr1}
and \eqref{eq:z0coeff} gives Theorem \ref{thm:main}.

\subsection*{Acknowledgement:} The author thanks to András Némethi for fruitful conversations about the problem. The author was partially supported by the \emph{Lend\"ulet program} of the Hungarian Academy of Sciences and by the ERC Advanced Grant LDTBud (awarded to Andr\'as Stipsicz).

\bibliographystyle{amsplain}
\bibliography{antidiag}

\end{document}